\documentclass[12pt]{article}

\title{A new proof of Seymour's 6-flow theorem}

\author{
   Matt DeVos\thanks{
     Email: {\tt mdevos@sfu.ca}. 
     Supported by an NSERC Discovery Grant (Canada).}
 \and
   Edita Rollov\'{a}\thanks{
  Email: {\tt rollova@ntis.zcu.cz}. 
  European Centre of Excellence, NTIS – New Technologies for Information Society, Faculty of Applied Sciences, University of West Bohemia, Pilsen. 
  Partially supported by project GA14-19503S of the Czech Science Foundation. 
  Partially supported by project LO1506 of the Czech Ministry of Education, Youth and Sports. 
}
 \and
   Robert \v{S}\'{a}mal\thanks{
  Email: {\tt samal@iuuk.mff.cuni.cz}. 
  Computer Science Institute of Charles University, Prague. 
  Partially supported by grant GA \v{C}R P202-12-G061.
  Partially supported by grant LL1201 ERC CZ of the Czech Ministry of Education, Youth and Sports.
}
}

\date{}

\usepackage{fullpage}
\usepackage{graphicx}
\usepackage{amsfonts}
\usepackage{amsmath}
\usepackage{amsthm}
\usepackage{enumitem}

\theoremstyle{plain}
\newtheorem{theorem}{Theorem}
\newtheorem{lemma}[theorem]{Lemma}

\newtheorem{observation}[theorem]{Observation}

\newtheorem{conjecture}[theorem]{Conjecture}

\theoremstyle{definition}
\newcommand\supp{\operatorname{supp}}

\begin{document}

\maketitle

\begin{abstract}
Tutte's famous 5-flow conjecture asserts that every bridgeless graph has a nowhere-zero 5-flow.  Seymour proved that every such graph has a nowhere-zero 6-flow.  Here we
give (two versions of) a new proof of Seymour's Theorem.
Both are roughly equal to Seymour's in terms of complexity, but they offer an alternative perspective which we hope will be of value.  
\end{abstract}

\section{Introduction}
By default we permit graphs to have parallel edges, and we forbid loops. Let $G$ be a graph which is equipped with an orientation of its edges and let $\Gamma$ be an
additively written abelian group.  We define a function $\varphi : E(G) \rightarrow \Gamma$ to be a $\Gamma$-\emph{flow} if the following condition, called
\emph{Kirchhoff's law}, is satisfied at every vertex $v$.  Here $\delta^+(v)$ ($\delta^-(v)$) denotes the set of edges directed away from (toward) $v$.
\[ \sum_{e \in \delta^+(v)} \varphi(e) = \sum_{e \in \delta^-(v)} \varphi(e) \]
We say that $\varphi$ is \emph{nowhere-zero} if $0 \not\in \varphi(E(G))$ and for a positive integer $k$ we say that $\varphi$ is a  $k$-\emph{flow} if $\Gamma = \mathbb{Z}$ and $| \varphi(e) | < k$ holds for every edge $e$.  Note that if $\varphi$ is a $\Gamma$-flow, we may reverse the orientation of an edge $e$ and replace $\varphi(e)$ by its additive inverse to obtain a new flow.  This operation preserves the properties of being nowhere-zero and $k$-flow, so the question of whether $G$ has a nowhere-zero $\Gamma$-flow ($k$-flow) depends only on the underlying graph and not the particular orientation.  Accordingly, we say that an undirected graph $G$ has a nowhere-zero $\Gamma$-flow ($k$-flow) if some (and thus every) orientation of $G$ admits such a function. 

Tutte initiated the study of nowhere-zero flows by demonstrating that nowhere-zero $k$-flows are dual to $k$-colourings in planar graphs.  He made three famous conjectures, known as the 5-flow, 4-flow, and 3-flow conjectures which extend various colouring theorems about planar graphs to arbitrary graphs.  These conjectures have been the driving motivation for the subject, and all three remain open despite considerable effort.  For the purposes of this paper, the relevant conjecture is the 5-flow conjecture.
The reader interested to learn more about the area may consult~\cite[Chapter~6]{Diestelbook} or~\cite{Zhangbook}. 

\begin{conjecture}[Tutte \cite{WT}]
Every 2-edge-connected graph has a nowhere-zero 5-flow.
\end{conjecture}

The Petersen graph does not have a nowhere-zero 4-flow, so if true, this conjecture is best possible.  The strongest partial result toward this conjecture is Seymour's 6-flow theorem.

\begin{theorem}[Seymour \cite{PS}]
\label{seymour}
Every 2-edge-connected graph has a nowhere-zero 6-flow.
\end{theorem}

The purpose of this paper is to provide a new proof of this theorem.  In fact we will give two (closely related) variants of this new proof.  Seymour's original paper also
contains two proofs of the 6-flow theorem which are similar in some ways, but have somewhat different character.  Both of our proofs are roughly equal to Seymour's in
terms of complexity, but they offer an alternative perspective which we hope will be of value.  All of these proofs rely on the following key result of Tutte.

\begin{theorem}[Tutte]
\label{int2mod}
For every graph $G$ and positive integer $k$ we have the following equivalence: $G$ has a nowhere-zero $k$-flow if and only if $G$ has a nowhere-zero $\mathbb{Z}_k$-flow.
\end{theorem}

Though we shall not require it, let us note that Tutte proved more generally that for every abelian group $\Gamma$ of order $k$, a graph has a nowhere-zero $k$-flow if and
only if it has a nowhere-zero $\Gamma$-flow.  For our purposes, the key consequence of Theorem \ref{int2mod} (and the isomorphism 
$\mathbb{Z}_6 \cong \mathbb{Z}_2 \times \mathbb{Z}_3$) is that the problem of finding a nowhere-zero 6-flow in a graph $G$ is now reduced to that of finding a
$\mathbb{Z}_2$-flow $\varphi_2$ and a $\mathbb{Z}_3$-flow $\varphi_3$ for which $\supp(\varphi_2) \cup \supp(\varphi_3) = E(G)$.

In addition to this reduction to modular flows, all present proofs of the 6-flow theorem rely upon some standard reductions to 3-connected cubic graphs; this is the focus of the following section.  For now, let us take a moment to compare and contrast these different proofs.  Seymour's first proof involves finding a sequence of cycles in the graph $C_1, \ldots, C_t$ with the property that for every $j > 1$ there are at least two edges with one end in $V(C_j)$ and one end in $\cup_{i=1}^{j-1} V(C_i)$.  These cycles are precisely the support of the $\mathbb{Z}_2$-flow.  So, we may view this proof as a process in which we repeatedly decide upon the values of the $\mathbb{Z}_2$-flow on some edges, and a certain graph structure is guaranteed to permit the construction of a suitable $\mathbb{Z}_3$-flow at the finish.  Our proofs do the opposite of this.  Namely, in our process we repeatedly decide upon the values of the $\mathbb{Z}_3$-flow on certain edges, and maintain a structure which permits us to find a suitable $\mathbb{Z}_2$-flow at the end.  

Seymour's second proof of his 6-flow theorem relies upon an argument somewhat similar to the first, but instead constructs two disjoint sets of edges with special
structure.  One is the edge-set of a spanning tree, the other is called a 2-base.  One may construct a $\mathbb{Z}_2$-flow whose support contains the edges in the
complement of a spanning tree, and a $\mathbb{Z}_3$-flow whose support contains the edges in the complement of a 2-base, and this gives the desired result.  In fact, this
proof shows that a more restrictive property, called $\mathbb{Z}_6$-\emph{connectivity}~\cite{JLPT} holds for 3-edge-connected graphs.  
Namely, for every oriented 3-edge-connected graph $G$ and every function $f : E(G) \rightarrow \mathbb{Z}_6$ there exists a flow $\varphi : E(G) \rightarrow \mathbb{Z}_6$
for which $\varphi(e) \neq f(e)$ holds for every $e \in E(G)$.  This group connectivity analogue of the 6-flow theorem can also be obtained by straightforward modification
of our proofs, though we will not provide these details.

\section{A Standard Reduction}

We will use standard techniques to reduce the proof of the 6-Flow Theorem to that of proving the following lemma.  (This is a key step in both of Seymour's proofs and also both of ours.)

\begin{lemma}[Seymour]
\label{reduction}
Every 3-edge-connected cubic graph has a nowhere-zero $\mathbb{Z}_6$-flow.
\end{lemma}

So, in this section we will prove Seymour's 6-Flow Theorem assuming the above lemma.  Then in each of the following sections we give a new proof of this lemma.  We begin with a
straightforward observation. (Here $\delta^+(X)$ ($\delta^-(X)$) denotes the set of edges with initial vertex in $X$ (in $V(G) \setminus X$) and terminal vertex 
in~$V(G) \setminus X$ (in $X$).)

\begin{observation}
\label{edge-cut}
Let $\Gamma$ be an abelian group, let $G=(V,E)$ be an oriented graph and let $X \subseteq V$.   If
$\varphi : E \rightarrow \Gamma$ satisfies Kirchhoff's law at every $x \in X$ then 
$\sum_{e \in \delta^+(X)} \varphi(e) = \sum_{e \in \delta^-(X)} \varphi(e)$. 
\end{observation}

\begin{proof} 
$ 0 = \sum_{x \in X} \left( \sum_{e \in \delta^+(x)} \varphi(e) - \sum_{e \in \delta^-(x)} \varphi(e) \right) = \sum_{e \in \delta^+(X)} \varphi(e) - \sum_{e \in \delta^-(X)} \varphi(e)$.
\end{proof}

Note that by this observation, whenever $\varphi : E(G) \rightarrow \Gamma$ satisfies Kirchhoff's law at every vertex except possibly $v$, then $\varphi$ also satisfies
Kirchhoff's law at $v$ (take $X = V(G) \setminus \{v\}$).  With this observation and Lemma \ref{reduction}, we are ready to prove the 6-Flow Theorem.%It follows from this that whenever we have a flow $\varphi$ on a graph $G$ and we modify $G$ by contracting the edge $e$, then the restriction of $\varphi$ to $E(G) \setminus \{e\}$ gives a flow in this new graph.  

\begin{proof}[Proof of Theorem \ref{seymour}]
By Theorem \ref{int2mod}, it suffices to show that every 2-edge-connected graph~$G$ has a nowhere-zero $\mathbb{Z}_6$-flow.  Suppose (for a contradiction) this is false
and choose a counterexample $G$ which is minimal in the following sense (here $V_{3+}(G) = \{ v \in V(G) \mid \deg(v) \ge 3 \}$): 
\begin{enumerate}[label=(\arabic*)]
\item The number of edges of $G$ in a 2-edge-cut is minimum,
\item $\sum_{v \in V_{3+}(G)} (\deg(v) - 3 )$ is minimum, subject to (1).
\end{enumerate}

We may assume $|V(G)| \ge 3$ as otherwise the result holds trivially.  
If $G$ has a cut-vertex $v$, then there exist nontrivial edge-disjoint subgraphs $G_1,G_2$ with $G_1 \cup G_2 = G$ and $V(G_1) \cap V(G_2) = \{v\}$.  Since $G_1$ and $G_2$ are 2-edge-connected, it follows from our minimality criteria that both $G_1$ and $G_2$ have nowhere-zero $\mathbb{Z}_6$-flows. But then the same is true for $G$, a contradiction.  

If $G$ has a 2-edge-cut $\{e,e'\}$, then by the minimality of our counterexample, $G/e$ has a nowhere-zero $\mathbb{Z}_6$-flow $\varphi$.  It follows from basic principles that we may assign a value $\varphi(e)$ so that $\varphi$ is a flow in $G$.  By Observation \ref{edge-cut}, $\varphi(e) = \pm \varphi(e')$, so $\varphi$ is nowhere-zero, a contradiction.

Next suppose that $G$ has a vertex~$v$ with $\deg(v) > 3$, and let $e_1, \ldots, e_k$ be the
edges incident with $v$.  Form a new graph $G'$ from $G$ by adding a disjoint cycle $C$ with $V(C) = \{ v_1, \ldots, v_k \}$, and changing the edge $e_i$ to have endpoint
$v_i$ instead of $v$.  After this modification $v$ is an isolated vertex and we delete it. 
Now, if $G'$ contains a 1-edge-cut or a 2-edge-cut, then it is easy to find in~$G$ a 1-edge-cut or a 2-edge-cut or a cut-vertex. 
Consequently, $G'$ is 3-edge-connected, so by minimality it has a nowhere-zero $\mathbb{Z}_6$-flow.  
Identifying $V(C)$ to a vertex  returns us to the original graph~$G$, where we have found a nowhere-zero $\mathbb{Z}_6$-flow, a contradiction.

We have shown that $G$ is cubic and 3-edge-connected, so the result follows from Lemma~\ref{reduction}.
\end{proof}

\section{First Proof}

In this section in some instances we will form a new graph from an old one by identifying a set $X$ of vertices to a single new vertex.  In order to stay within the realm of loopless graphs, we shall always assume that any edge with both ends in $X$ will be deleted by this process.  

Let $G$ be a 2-edge-connected graph, and define a relation $\sim$ on $E(G)$ by the rule $e \sim f$ if either $e = f$ or $\{e,f\}$ is a 2-edge-cut.  It follows from
elementary principles that this is an equivalence relation, and we call its equivalence classes \emph{generalized series classes}.  It is an easy exercise to show that for
a generalized series class $F$ with $|F| \ge 2$, the graph obtained from $G$ by identifying each component of $G-F$ to a single vertex is a cycle
(possibly of length~2). The following lemma is the main result of this section. (Here $\delta(u)$ denotes the set of edges incident to the vertex $u$.)

\begin{lemma}
Let $G=(V,E)$ be an oriented graph with a distinguished root vertex $u$, let $S \subseteq \delta(u)$ satisfy $|S| = 2$ and let $\psi_2 : S \rightarrow \mathbb{Z}_2$ and
$\psi_3 : \delta(u) \rightarrow \mathbb{Z}_3$ be functions.  Suppose 
\begin{enumerate}[label=(\roman*)]
\item $\deg(v) = 3$ for every $v \in V \setminus \{u\}$, 
\item $G$ is 3-edge-connected, 
\item $G-u$ is 2-edge-connected, 
\item $\sum_{e \in \delta^+(u)} \psi_3(e) - \sum_{e \in \delta^-(u)} \psi_3(e) = 0$.
\end{enumerate}
Then there exist flows $\varphi_2 : E \rightarrow \mathbb{Z}_2$ and $\varphi_3 : E \rightarrow \mathbb{Z}_3$ so that 
$\varphi_2|_S = \psi_2$ and $\varphi_3|_{\delta(u)} = \psi_3$ and $(\varphi_2(e), \varphi_3(e)) \neq (0,0)$ for every $e \in E \setminus \delta(u)$.  
\end{lemma}

\begin{proof}
We proceed by induction on $|V|$.  For the base case when $|V| = 2$, we set $\varphi_3 = \psi_3$ and we choose $\varphi_2 : E \rightarrow \mathbb{Z}_2$ so that
$\varphi_2|_S = \psi_2$ and $|\supp(\varphi_2)|$ is even (this is always possible since the first condition guarantees $\delta(u) \setminus S \neq \emptyset$).  It follows from basic principles that $\varphi_2$ and $\varphi_3$ are flows and the remaining condition is trivially satisfied.

For the inductive step, begin by choosing $uv \in \delta(u) \setminus S$ and define $F$ to be the generalized series class in $G-u$ which contains both edges incident with 
$v$. (Note that this is always possible since $u$ is of degree at least three by (ii).) Let $G_0, \ldots, G_{\ell-1}$ be the components of $(G-u) - F$ and assume that they are ordered so that $F$ contains an edge, called $f_k$, between $G_k$ and 
$G_{k+1}$ (working mod $\ell$) for every $0 \le k \le \ell-1$. The edge~$f_k$ is defined uniquely, except for the case $l=2$, when there are two 
edges between $G_0$ and~$G_1$. In this case we require that $f_0 \ne f_1$. 
Since $G$ is 3-edge-connected, there is an edge from~$u$ to each of the graphs~$G_k$. 

Assume further that $G_0$ consists of the single vertex $v$ and that one edge in $S$ has an end in $G_i$ 
and the other edge in $S$ has an end in $G_j$ where $0 < i \le j$.  

\begin{figure}[h]
  \centerline{\includegraphics{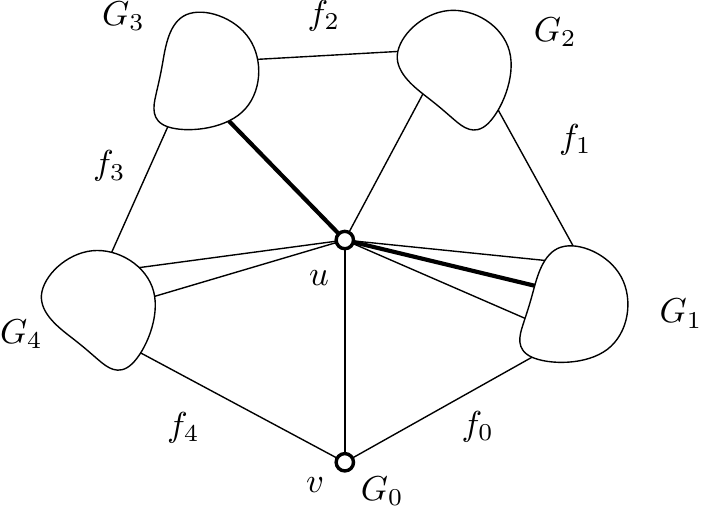}}
  \caption{Illustration of our first proof. We have $l=5$, $i=1$, $j=3$, the edges in~$S$ are drawn in bold.} 
\end{figure}

Next we shall extend the domains of our functions $\psi_3$ and $\psi_2$.  Let $G^*$ be the graph obtained from $G$ by identifying $V(G_i)$ to a single vertex for every $0
\le i \le \ell-1$.  It follows from basic principles that we can extend the domain of $\psi_3$ to $F$ in such a way that $\psi_3$ is a $\mathbb{Z}_3$-flow on $G^*$.  Moreover, there is a $\mathbb{Z}_3$-flow of $G^*$ with support $F$, and by adding a suitable multiple of this flow to $\psi_3$, we may further assume that $\psi_3(f_{i-1}) \neq 0 \neq
\psi_3(f_j)$.  Extend the domain of $\psi_2$ by defining $\psi_2(f) = 1$ for all $f \in F \setminus \{ f_{i-1}, f_j \}$.  (We  postpone deciding $\psi_2(f_{i-1})$ and
$\psi_2(f_j)$ for the moment.)

For every $0 \le k \le \ell-1$, form the graph $G_k^+$ from the original graph $G$ by identifying all vertices in $V \setminus V(G_k)$ to a single root vertex.  For each
$G_k^+$, the restriction of $\psi_3$ to the edges incident with the root vertex obeys Kirchhoff's law at this vertex.  To see this, note that this is equivalent to the
statement that this restriction is a flow in the two vertex graph obtained from $G_k^+$ by identifying $V(G_k)$ to a single point.  However, this latter graph may also be
obtained by identifying vertices in $G^*$, where $\psi_3$ is a flow.  

If $i=j$ then both edges in $S$ have an endpoint in $G_i$ and we apply the lemma inductively to the graph $G_i^+$ together with the set $S$ and the functions $\psi_2|_S$
and $\psi_3$ restricted to the edges incident with the root of $G_i^+$.  Let $\varphi_2^i : E(G_i^+) \rightarrow \mathbb{Z}_2$ and $\varphi_3^i : E(G_i^+) \rightarrow
\mathbb{Z}_3$ be the resulting flows.  Define $\psi_2(f_{i-1}) = \varphi_2^i(f_{i-1})$ and $\psi_2(f_{i}) = \varphi_2^i(f_i)$.  Now for every other graph $G_k^+$ with $k
\neq i$ we apply the lemma inductively with the set $\{f_{k-1}, f_k \}$ in place of~$S$ and the appropriate restrictions of $\psi_2$ and $\psi_3$.  The $\mathbb{Z}_2$-flows and
$\mathbb{Z}_3$-flows obtained by these inductive applications of the lemma combine to yield a solution for the original graph:
clearly, they extend $\psi_2$ and~$\psi_3$, and satisfy the nonzero condition. They satisfy Kirchhoff's law at every vertex of some~$G_k^+$, thus 
at every vertex of~$V \setminus \{u \}$. Therefore, they satisfy Kirchhoff's law everywhere. 

If $i < j$ then we apply induction to the graph $G_i^+$ together with the two-edge set consisting of the edge $f_i$ and the edge in $S$ incident with $G_i$ (together with
the restrictions of $\psi_3$ and $\psi_2$ to these sets).  Let $\varphi_2^i : E(G_i^+) \rightarrow \mathbb{Z}_2$ and $\varphi_3^i : E(G_i^+) \rightarrow \mathbb{Z}_3$ be
the resulting flows.  Similarly, we apply induction to the graph $G_j^+$ together with the two edge set consisting of $f_{j-1}$ and the edge in $S$ incident with $G_j$ and
we let $\varphi_2^j : E(G_j^+) \rightarrow \mathbb{Z}_2$ and $\varphi_3^j : E(G_j^+) \rightarrow \mathbb{Z}_3$ be the resulting flows.  Define $\psi_2(f_{i-1})
= \varphi_2^i(f_{i-1})$ and $\psi_2(f_j) = \varphi_2^j(f_j)$.  (Here we use the fact that $0 < i < j$, thus $f_{i-1}$ and~$f_j$ are distinct edges.) 
As before, for every other graph $G_k^+$ with $k \not\in \{i,j\}$ we apply the lemma inductively with the
set $\{f_{k-1}, f_k \}$ in place of~$S$ and the appropriate restrictions of $\psi_2$ and $\psi_3$.  The $\mathbb{Z}_2$-flows and $\mathbb{Z}_3$-flows obtained by these
inductive applications of the lemma combine to yield a solution for the original graph. 
\end{proof}

\begin{proof}[Proof of Lemma \ref{reduction}]
Let $G$ be a 3-edge-connected cubic graph and apply the above lemma to $G$ with an arbitrary root vertex~$u$ and functions $\psi_2, \psi_3$ for which 
$\supp(\psi_2) \cup \supp(\psi_3) = \delta(u)$. (Note that $G-u$ is 2-edge-connected, since $\deg(u)=3$.) 
\end{proof}

\section{Second Proof}
For our second proof we will work with more general functions than flows and for this we will require some additional notation.  For a function $\varphi : E \rightarrow
\mathbb{Z}_k$, the \emph{boundary} of $\varphi$ is the function $\partial \varphi : V \rightarrow \mathbb{Z}_k$ given by the rule
\[ \partial \varphi (v) = \sum_{e \in \delta^+(v)} \varphi(e) - \sum_{e \in \delta^-(v)} \varphi(e). \]
Every such function satisfies the \emph{zero-sum rule}: $\sum_{v \in V} \partial \varphi(v) = 0$ since (after expanding) each edge $e$ contributes 
$\varphi(e) - \varphi(e) = 0$ to the total.  

Working with general functions instead of flows has the advantage that we may use a straightforward induction step: delete an edge $vw$ and modify 
boundaries of~$v$ and of~$w$. This is similar in spirit to recent proof of Thomassen et al.~\cite{Thom}, albeit much easier. 
The precise statement we will be proving is given by the following lemma. 
(Here $V_t(G)$ denotes the set of vertices of degree~$t$ in the graph~$G$, and we call $G$ \emph{subcubic} if 
$V_t(G) = 0$ for all $t > 3$.)

\begin{lemma}
Let $G$ be an oriented 2-edge-connected subcubic graph with root vertex $u \in V_2(G)$.  Let $\mu : V(G) \rightarrow \mathbb{Z}_3$ and for $k=2,3$ let $\psi_k : \delta(u) \rightarrow \mathbb{Z}_k$.  Suppose 
\begin{enumerate}[label=(\roman*)]
\item $\sum_{v \in V} \mu(v) = 0$, 
\item $\supp(\mu) \subseteq V_2(G)$, 
\item $\partial \psi_3(u) = \mu(u)$,
\item $\partial \psi_2(u) = 0$ if $\mu(u) = 0$, and
\item $( \psi_2(e), \psi_3(e) ) \neq (0,0)$ for every $e \in \delta(u)$.  
\end{enumerate}
Then there exist functions $\varphi_k : E \rightarrow \mathbb{Z}_k$ for $k=2, 3$ satisfying 
\begin{enumerate}[label=(\arabic*)]
\item $\varphi_k|_{\delta(u)} = \psi_k$ for $k=2,3$, 
\item $\partial \varphi_3 = \mu$, 
%\item $\supp(\partial \varphi_2) \subseteq \supp(\mu)$,
\item for every vertex~$v$ holds that $\partial \psi_2(v) = 0$ if $\mu(v) = 0$, and, finally, 
\item $(\varphi_2(e), \varphi_3(e)) \neq (0,0)$ for every $e \in E(G)$.
\end{enumerate}
\end{lemma}

\begin{proof}
We proceed by induction on $|E(G)|$.  Our first base case is when $G$ is a cycle of length~2.  Let $V(G) = \{u,v\}$ and note that 
$\partial \psi_k(v) = - \partial \psi_k(u)$ by the zero-sum rule.  
Now (i) with (iii) imply~(2) % that $\partial \psi_3 = \mu$ 
and (i) with (iv) imply~(3). %that $\supp(\partial \psi_2) \subseteq \supp(\mu)$.  
So $\varphi_2 = \psi_2$ and $\varphi_3 = \psi_3$ satisfy the lemma.  We may now assume that there is no vertex $v \in V_2(G) \setminus  \{u\}$ with
$\mu(v) = 0$, since otherwise the result follows by applying induction to a graph obtained by contracting an edge incident with $v$ but not $u$.  
Our second base case is when $G$ is a cycle of length 3.  Let $V(G) = \{u,v,w\}$ and note that by assumption $\mu(v) \neq 0 \neq \mu(w)$.  
For $k=2,3$ extend $\psi_k$ to a function $\varphi_k$ defined on $E(G)$ as follows.   
For $k=3$, choose $\varphi_3(vw)$ so that $\partial \varphi_3 = \mu$ (this is possible by (i), (iii) and the zero-sum rule), 
for $k=2$ set $\varphi_2(vw) = 1$.  The resulting functions satisfy the lemma.  

Next suppose that $G$ has an edge-cut $S$ of size two which separates the vertices into $X_1,X_2$ where $|X_1|, |X_2| \ge 2$.  Assume (without loss) that $u \in X_2$ and for $i=1,2$ let $G_i$ be the graph obtained from $G$ by identifying $X_i$ to a new vertex $x_i$ (and deleting any resulting loops).  For $i=1,2$ define $\mu^i : V(G_i) \rightarrow \mathbb{Z}_3$ by the rule
\[ \mu^i(v) = \left\{ \begin{array}{cl}
	\mu(v)	&	\mbox{if $v \neq x_i$}	\\
	\sum_{x \in X_i} \mu(x)	&	\mbox{if $v = x_i$.} 
	\end{array} \right. \]
Apply induction to $G_1$ together with $\mu^1$, the root vertex $u$ and $\psi_2$, $\psi_3$ to obtain $\varphi^1_2, \varphi^1_3$.  
Next apply induction to $G_2$ with $\mu^2$, the root vertex $x_2$ and the functions $\varphi^1_2|_S$ and $\varphi^1_3|_S$. 
Observe that $\mu^2(x_2) = - \mu^1(x_1)$ and thus (3)~for~$G_1$ implies (iv)~for~$G_2$. 
The second application of induction gives us functions~$\varphi^2_2$ and~$\varphi^2_3$. 
Merging the functions~$\varphi^1_k$ and~$\varphi^2_k$ for $k=2,3$ yields a solution.  
So we may now assume that every 2-edge-cut of $G$~is trivial (i.e., associated with a single vertex).

Next suppose there exists a vertex $v \in V_2(G) \setminus  \{u\}$ and let $w_1,w_2$ be its neighbours (without loss of generality we assume that the edges $vw_1$ and $vw_2$ are oriented away from~$v$).  
Choose  $\tau : \{ vw_1, vw_2 \} \rightarrow \mathbb{Z}_3 \setminus \{0\}$ so that $\partial \tau (v) = \mu(v)$.  Then define $\mu'  : V(G-v)  \rightarrow \mathbb{Z}_3$ by
the rule $\mu'(w_i) = \mu(w_i) - \partial \tau (w_i)$ for $i=1,2$ and otherwise $\mu'(w) = \mu(w)$.  It follows from~(i) and the zero-sum rule for $\tau$ that $\sum_{w \in
V \setminus \{v\}} \mu'(w) = 0$.  Moreover, since every 2-edge-cut of $G$ is trivial, $G-v$ is 2-edge-connected. So, we may apply induction to $G - v$ together with $\mu'$, $\psi_2$, and $\psi_3$ to obtain $\varphi'_2$ and $\varphi'_3$.  Extend
$\varphi'_3$ to a function $\varphi_3 : E(G) \rightarrow \mathbb{Z}_3$ by defining $\varphi_3(vw_i) = \tau(vw_i)$ for $i=1,2$.  Extend $\varphi'_2$ to a function
$\varphi_2 : E(G) \rightarrow \mathbb{Z}_2$ by defining $\varphi_2( vw_i ) = \partial \varphi'_2 (w_i)$.  Now $\varphi_2, \varphi_3$ give a solution.  

We may now assume that $V_2(G) = \{u\}$ and $\mu = 0$.  Choose an edge $vw$ with $v,w \neq u$, and note that $G-vw$ is 2-edge-connected.  Define $\mu' : V(G) \rightarrow \mathbb{Z}_3$ by 
$\mu'(v) = 1$, $\mu'(w) = -1$ and $\mu'(x) = 0$ for all $x \in V(G) \setminus \{v,w\}$.  
Apply induction to $G-vw$ together with $\mu'$, $\psi_2$, and $\psi_3$ to obtain functions $\varphi_2$ and $\varphi_3$.  
Now we may extend the domains of these functions by choosing $\varphi_3(vw) = \pm 1$ and $\varphi_2(vw) \in \mathbb{Z}_2$ so that $\partial \varphi_k = 0$ for $k=2,3$.
This yields the desired functions for $G$.
\end{proof}

\begin{proof}[Proof of Lemma \ref{reduction}]
Let $G=(V,E)$ be an oriented 3-edge-connected cubic graph, and modify $G$ by subdividing an edge with a new root vertex $u$.  
Let $\mu$ be identically $0$ and choose $\psi_k \colon  \delta(u) \rightarrow \mathbb{Z}_k$ for $k=2,3$ so that $\partial \psi_k (u) = 0$ and $( \psi_2(e), \psi_3(e) ) \neq (0,0)$ for every $e \in \delta(u)$.  Now applying the lemma (and suppressing $u$) gives a nowhere-zero $\mathbb{Z}_6$-flow of $G$, as desired.
\end{proof}

\end{document}